# Optimal deployment of collecting routes for urban recycle stations (eco-points).

G. Marseglia[ab*], J.A. Mesa[bcd], F.A. Ortega[e], R. Piedra-de-la-Cuadra[e].

a) Research Department, Link Campus University of Rome, Via del Casale di San Pio V, 44-00165, Rome Italy

b) Universidad de Sevilla, (IMUS), Universidad de Sevilla, Avenida Reina Mercedes 41012 Seville, Spain

c) Higher Technical School of Engineering, Universidad de Sevilla, C/ Américo Vespucio 41092 Seville, Spain.

d) Departamento de Matemática Aplicada II, Universidad de Sevilla, Camino de los Descubrimientos s/n, 41092, Seville, Spain.

d) Higher Technical School of Architecture, Universidad de Sevilla, Avenida Reina Mercedes 41012 Seville, Spain.

*Corresponding author e-mail: g.marseglia@unilink.it

## Abstract

The rapid and constant increase in urban population has led to implied a drastic rise in urban solid waste production with worrying consequences for the environment and society. In many cities, an efficient waste management combined with a suitable design of vehicle routes (VR) can lead to benefits in the environmental, economic, and social impacts.

The general population are becoming increasingly aware of the need for the separation of the various categories of municipal solid waste. The numerous materials collected include glass, PET or batteries, and electric components, which are sorted at the eco-points. The management of eco-points gives rise to several problems that can be formulated analytically. The location and number of eco-point containers, the determination of the fleet size for picking up the collected waste, and the design of itineraries are all intertwined, but present computationally difficult problems, and therefore must be solved in a sequential way.

A mathematical model has been formulated in this paper, based on the Bin Packing (BP) and VR schemes, for the deployment of routes of mobile containers in the selective collection of urban solid waste. A heuristic algorithm has also been developed, which considers two different configurations of the containers to solve the proposed mathematical programming model. The results obtained from the numerical simulations show the validation of the proposed methodology carried out for the Sioux Falls network benchmark and the specific real case study.

## 1. Introduction

Nowadays, due to greater attention to the quality of life and the interest in sustainable energy resource utilisation, the majority of institutions are moving towards efficient solutions to combat the urban waste collection issue. The Sustainable Development Goals defined by the United Nations in the 2030 Agenda have brought the imperative to make cities more sustainable to the fore of development discussions, thereby improving the quality of life for citizens [1]. European Nations have recently defined several policies in favour of the circular economy [2,3]. Therefore, due to the inherent social implications, the optimisation of urban waste collection assumes a fundamental role in each city today [4].

In recent years, the growth in urban population density has implied a major rise in the production of various kinds of Municipal Solid Waste (MSW), whose management includes several functional phases, such as waste generation, storage, collection, transportation, processing, recycling, and disposal in a suitable landfill. As a consequence, administrations, such as municipalities, have defined suitable waste collection areas to obtain efficiency and low environmental impact. On the other hand, the detailed analysis of the waste management strategy can also lead to major benefits in terms of cost reduction. In fact, urban waste allocation and transportation influence the costs of MSW by approximately 60%-80% [5,6]. Agovino et al. [7] analysed the case of Italian separate waste collection. They highlighted the evolution of recent European waste policies to environmental action plans and legislation framework, for the attainment of environmental benefits and of an energy- and resource-efficient economy.

In addressing this strong link between this type of issue and the benefits for citizens, there is a growing awareness regarding the definition of the best management and logistics system. Solid-Waste Management (SWM) is focused on providing solid-waste collection services for residential, industrial, and commercial customers in metropolitan areas. Residential customers are generally people living in private homes, whereas commercial customers include strip malls, restaurants, and small office buildings. The choice of permanent or temporary collection points in the various areas of each city is of extreme importance. Indeed, every waste container presents a specific capacity, cost, and environmental impact. In this respect, the actual challenge for most of the administrations involves the suitable definition of key parameters, such as the type, the number, and the position of the containers for every area in order to dispose of all types of waste produced in a defined period. A further important aspect for consideration in the analysis of this issue is that of the size of the city examined. In fact, in highly populated cities, waste collection is managed by different municipalities or organisations responsible for specific established zones.

Therefore, the constant improvement in waste recovery in terms of demand and technology factors and in the mentality of citizens is addressed by an increasing number of researchers for the identification of the optimal management design capable of meeting future requirements. In this scenario regarding municipalities, it will also be problematical to predict future modifications in waste collection and recycling. The composition of municipal solid urban waste is influenced by the standard of living of the population, the economic activity of their inhabitants, and the climate of the region [8]. Certain products will eventually become more commonly used in relation to these factors and, subsequently, various waste modalities will be generated in varying proportions.

According to the report entitled Statistics on Waste Collection and Treatment for the Year 2018 of the National Institute of Statistics (INE) in Spain, the main materials produced in Spain are paper and cardboard (24.1%), organic matter (22.9%), glass (18.9%), plastic and mixed packaging (16.8%), and others (representing 17.3%). This latter group requires special attention, since certain items can be considered as hazardous waste. Waste can be characterised as hazardous if it possesses any one of the following four features: ignitability, corrosiveness, reactivity, and toxicity. Hazardous waste, which is usually the waste by-product of our industrial processes, presents immediate or long-term risks to humans, animals, plants, or the environment.

In Spain, many municipalities have combined the need to collect this type of potentially hazardous waste with the promotion of environmental policies and use containers with an aesthetically attractive design, which help spread the commitment to the selective collection of urban solid waste.

The so-called eco-points are large waste containers with separate non-homogeneous sections for the collection of various kinds of items, including mobiles, batteries, chargers, syringes and needles, used low-energy lamps, radiographs and photographic material, books for the exchange between citizens, toner and cartridges of ink, aluminium and plastic coffee capsules, and CDs and DVDs A real example of an eco-point located in the city of Seville (Spain) is shown in Figure 1.

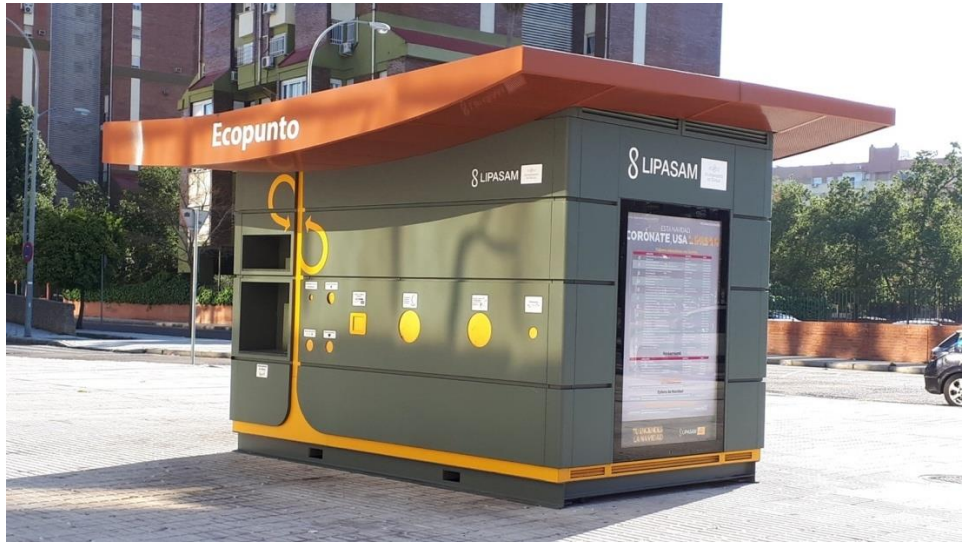

**Fig. 1**. An eco-point situated in city of Seville.

The decision regarding the best configuration for an eco-point container includes the design of the distribution of volume of their sections. The eco-point lay-out is linked to the Bin Packing (BP) problem. In agreement with Garey and Johnson [9], the BP issue is a combinatorial problem that belongs to the class of NP-hard problems. Various real applications of this kind of issue are presented in [10] and [11].

In addition to this connection with the BP problem, MSW collection is intrinsically connected to the VR model [12,13] in terms of optimising different criteria, such as the total distance travelled by vehicles, the emission of environmental pollutants, and the investment costs [14]. In this work, an optimisation model, inspired by the BP and VR schemes, has been formulated for the deployment of routes for mobile eco-points for the selective collection of solid municipal waste.

Our interest in this issue is focused on a type of mobile container for the collection of solid waste, composed of several sections for the separate storage of different items, which can either all be of the same size or can have heterogeneous volumes, depending on the demand of the place where they are temporarily located. Eco-points deployed in the region of Cantabria (Spain), such as that shown in Figure 2, represent real instances of the type of containers of interest. The identification of the allocation of the multi-block container is one of the two decisions to be adopted.

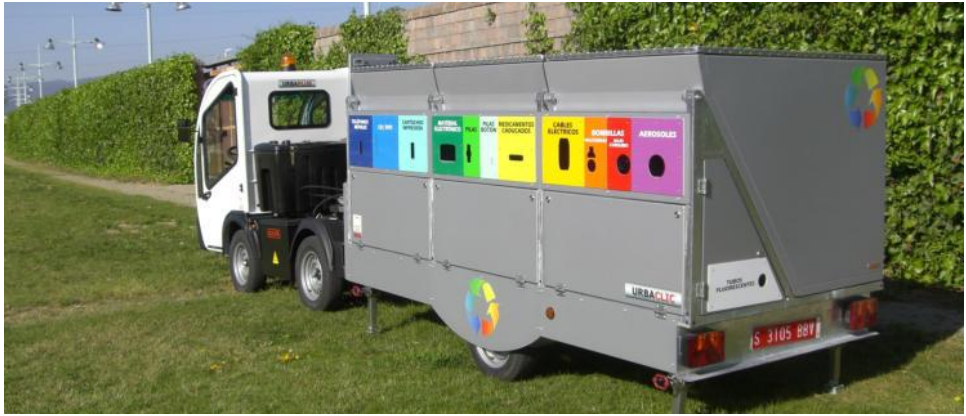

**Fig. 2**. An example of an eco-point deployed in Cantabria.

Beliën et al. [15] presented a review of the available literature on solid waste management problems where vehicle routing problems are classified into several categories. It should be borne in mind that mobile eco-points follow an established route and visit all the neighbourhoods of the city in an itinerant way. The container is placed at identified stops on public roads for a temporary period (for example, on a Monday, whereby it is moved to a new point in the city on the following Monday). The calendar is previously made known to all residents in the neighbourhood.

As stated in [17], the problem of managing the selective collection of waste containers within historic centres can be performed in three sequential phases: firstly, the location of containers along the streets; secondly, the determination of the minimum fleet size required to perform all collection services; and finally, the development of a model capable of defining the optimal paths, in terms of total and equilibrated number of kilometres travelled by the vehicles. As underlined in Barrena et al. [16], the result of the first phase (location of the containers) greatly influences the procedure since this determines the decision to be taken for the subsequent phases (route of collection vehicles and service programming). Recently, Barrena et al. [17] have developed a methodology for locating waste container facilities, whereby customer solidarity behaviour at each location is assumed, in order to produce results with a higher efficiency. Their model assumes that a proportion of inhabitants associated with a particular location node could take their waste to the closest container and another proportion of the population of the same node would be willing to take their waste to another unfilled container j' that is not excessively far away. This "solidary" behaviour of the clients would enable the efficient deployment of the containers in the area under analysis, thereby reducing their total number and grouping them at the points of lowest cost.

The main contribution of this work focuses on the double determination of BP configurations and of container routes. One of the ultimate aims involves the cost minimisation of the resources employed. In order to solve the problem of waste collection and vehicle routes in an optimal way, an adaptive algorithm of overflow deviated to the immediate neighbourhood is developed. This algorithm strives to solve the proposed mathematical programming model, whose computational complexity justifies the use of heuristics to address large real-life scenarios. The evaluation of the performance of the developed methodology has been carried out through a computational experience in a toy network using two strategies to design the layouts of the mobile multi-block containers that visit the demand nodes.

The remaining sections of the paper are organised as follows. In Section 2, previous research in this field is presented. The optimisation model, inspired in the BP and VR schemas, has been formulated for the deployment of routes of mobile eco-points for the selective collection of urban solid waste (Section 3). Subsequently, the formulation of this decision model has been developed, in addition to its corresponding resolution algorithm in Section 4. The performance evaluation of the proposed methodology based on a computational experience Sioux Falls network benchmark in a neighbourhood of Seville is highlighted in Section 5. Finally, several conclusions are summarised in Section 6.

## 2. Scientific Relevance and Literature Review

In this section, several papers are presented that are related to the issues of the Vehicle Route Problem (VRP) and of Bin Packing in waste collection management.

Waste collection and transportation problems are one of the most difficult operational issues in the development of an integrated waste management system [18,19]. Eiselt and Marianov [20] provide a survey of 64 studies on the landfill siting problem that include applications across the world, for which the main aspects of interest of the contributions are summarised in a table (see Ref. 20), by selecting country, technique, criteria, objectives, and type or facility to be located.

Furthermore, the intrinsic nature of MSW collection relates to the development of effective vehicle route models that optimise the total travelling distance of vehicles, the environmental emissions, and the investment costs [21]. Vehicle Routing is a scheduled process that enables vehicles to load waste at collection sites and to dump it at a landfill with the result being oriented by a single or multiple objectives [14].

In real-life scenarios, the waste collection system is distributed across a set of zones. Each zone has an associated starting and ending node. These nodes are used for vehicle routes and landfill points where the rubbish collected in the visited containers can be delivered. A planning horizon must also be considered in order to schedule a sequence of services within its bounds. A succession of routes (one per day, belonging to the same or to different distribution zones and performed by the same vehicle along the planning horizon) is called a circulation. Plans for the design of the vehicle circulation in railway transportation networks is described in [22].

Yeomans [23] underlined the importance of using a mathematical approach in the sustainable waste collection and transportation in urban areas. Melika Mohsenizadeh et al. [24] developed a mathematical model in order to reduce pollutant emissions from vehicles during the collection of waste, while considering the Ankara case. Beliën et al. [15] presented a review of solid waste management problems where VRPs are classified into several categories. Han and Ponce-Cueto [25] provided a detailed analysis of the Waste Collection Vehicle Routing Problem (WCVRP). The most efficient methods of VR and BP problems are based on heuristic and metaheuristic solution models [26-28]. Buhrkal et al. [29] applied an Adaptive Large Neighbourhood Search (ALNS) metaheuristic to solve the Waste Collection Vehicle Routing Problem with Time Windows (WCVRPTW) for various real cases. They obtained a minimisation of the distance driven which was also linked to a reduction in fuel consumption and pollutant emissions. Kim et al. [30] considered multiple disposal paths and drivers' lunch breaks using an extension of Solomon's approach [31]. Teixtera et al. [32] analysed the VRP while separating the collection of paper, glass, plastic, and metal waste materials. Their model is based on three stages: the definition of a zone for every truck; the identification of the type of waste to be collected every day; and the choice of the points to visit and of the sequence order. Angelelli and Speranza [33] presented a tabu search algorithm to optimise the operating costs of various waste collection processes for two case studies.

Bing et al. [34] investigated the plastic waste collection problem based on eco-efficiency in terms of the proper balance between environmental impacts, social issues, and costs. They modelled the urban plastic waste collection based on a VRP. Different cases were considered by analysing key factors, such as the type of truck, collection frequency, and collection method. These authors solved the VR issue by means of a tabu search algorithm capable of solving real cases. Their results showed that, based on the proposed algorithm, the eco-efficiency performance of the current collection paths could be improved by 7%. Das and Bhattacharyya [35] minimised the length of every waste collection and transportation route. They proposed a heuristic methodology based on the Travelling Salesman Problem (TSP) and obtained a reduction of approximately 30% of the total length of the waste collection route. Laporte et al. [36]

validated a heuristic algorithm for a Capacitated Arc Routing Problem (CARP) based on stochastic demand. More recently, Babaee Tirkolaee et al. [37] investigated the periodic CARP by considering the work of the drivers and their crews in order to analyse the demand change. Their model is based on an objective function for the minimisation of the total transport route and total costs in terms of the number of vehicles needed. A simulated annealing algorithm is used to improve the solution data.

A numerical approach based on Mixed Integer Linear Programming (MILP) is proposed in [38] and considers the VRP for every period to reduce the number of trucks, and in turn the total costs during a defined period. Mes et al. [39] presented a heuristic methodology based on the optimisation of various tunable parameters for every day of the week. Thanks to their model, they obtained a cost reduction of approximately 40% for a specific case study in the Netherlands. Son and Louati [40] developed a mathematical model that considered multiple transfer stations for urban solid waste management. They validated their model by applying it to a case study and obtained a reduction of path length and working hours of the trucks. Akhtar et al. [41] developed a modified Backtracking Search Algorithm (BSA) in Capacitated Vehicle Routing Problem (CVRP) models based on the "intelligent bin" idea to optimise path design in the waste management system. Their results for four days presented a 36.80% reduction in distance for 91.40% of the total waste collection, with an increase of the mean waste collection efficiency of 36.78% and a reduction in fuel consumption of 50%, in fuel cost of 47.77%, and in $CO_2$ production of 44.68%, respectively. Ghiani et al. [42] investigated the issue of locating waste collection points in urban areas. They developed a model to identify: the optimal sites for the location of the waste collection bins; the required number of bins; and the features of the bins sited at the various collection stations. The results obtained from their numerical model based on heuristic procedures demonstrated a reduction of approximately 62% for the waste collection points and a decrease of approximately 73% of the number of bins allocated. By considering capacity, time, and distance restrictions, Willemse and Joubert [26] extended the VRP in urban WM. They proposed four heuristic models to compute feasible solutions for the two issues of WM and VRP to minimise the total cost and/or the fleet dimension. They analysed CARP under time restrictions with intermediate facilities, and carefully modelled the collection of waste, based on the development of a constructive heuristic model.

Recently, various researchers have underlined the importance of optimising not only the logistic aspects of VRP, but also the bin packing methodologies to improve the performance of waste collection [43]. Rodrigues et al. [44] based the kerbside collection of WM on the classification of the bin components, truck components, and collection methods. They analysed the collection frequency by considering one bin for the single packaging waste flow (yellow for lightweight packaging; green for glass; blue for paper and cardboard). Martinho et al. [45] studied the importance of identifying the proper recycling methods, and analysed two different waste collection methods for the recycling in 3 districts of Portugal. In their research, they highlighted the importance of carrying out an itemised characterisation of the end waste to identify the overall quantity of waste packaging material that can be recycled. They defined certain both key "performance recycling indicators", such as waste characterisation, recycling rate, and separate waste collection rate, and also "logistic performance indicators", such as the distance travelled by car to collect waste, the number of workers involved in the collection, and the effective collection time.

The goal of the present work is not only to propose an optimisation model to cover the variable demand points in WM, but also to design the best VR in terms of minimal costs and environmental impact in the most efficient way. In this paper, the authors propose an algorithm for the optimisation of urban waste collection by means of routes for mobile eco-points for a random real example.

## 3 Model development

In order to determine optimal routes for mobile multi-block containers, we assume a strongly connected graph $G = (V, A)$, composed of a node set $V$ and an arc set $A$ (directed edges representing street sections), so that the existence of a shortest path in terms of distance (or travel time) between each pair of points of $V$ is guaranteed inside $G$. Let us suppose that set $V$ contains the set $J$ of nodes where waste is placed by the users in order to be collected in mobile multi-blocks containers. Each container is divided into sections (bins or blocks) and has associated one route, to be determined by the optimization model, that starts and ends at the same depot point ($O$) whose location is assumed to be fixed. Hence, the same index can be used to simultaneously represent each container and its corresponding route.

The following notation is used in our formulation:

    $I$:    Set of routes for containers ($i \in I$). We assume that all vehicles are identical and have the same transport capacity, that corresponds to one container. All containers are homogeneous and contain $|L|$ blocks of capacity $c$. Let $C$ be the total capacity of each container ($C = c \cdot |L|$).

    $J$:    Set of locations to be visited inside the city ($j, j' \in J \subseteq V$). Depot point $O$ is assumed to belong to this set $J$. Distances across network G between nodes $j$ and $j'$ are known and

recorded in $D = (d_{jj'})$ matrix. Term $d_{jj'}$ indicates the minimum cost of travelling from point $j$ to point $j'$. Once all the shortest paths between pairs of nodes of the set $J$ have been established, we can use, as a new solution space, a graph whose set of vertices is $J$ and where each arc $(j, j')$ is weighted by $d_{jj'}$. Let $A(J)$ be the set of these direct arcs between pairs of points in $J$. Analogously, we will denote by $A(S)$ the set of arcs that connect pairs of points belonging to subset $S \subseteq V$.

$K$: Set of waste modalities ($k \in K$).

Quantity of waste modality $k$ produced at point $j$ is represented by means of parameter $w_j^k \geq 0$.

Therefore:

$\sum_{k \in K} w_j^k$ : Indicates the total waste generated at point $j \in J$.

$\sum_{j \in J} w_j^k$ : Indicates the total waste of modality $k \in K$ produced inside the city.

$\sum_{k \in K} \sum_{j \in J} w_j^k$ Evaluates the total amount of waste produced inside the city for all waste modalities.

Moreover, the following variables are required in the model:

$y_i^k$ : Binary variable that takes value 1 if route $i \in I$ is used for collecting waste of modality $k \in K$, and is equal to 0 otherwise.

$x_{ij}^k$ : Binary variable that takes value 1 if location $j \in J$ is visited by route $i \in I$ and waste of modality $k \in K$ is collected, and is set to 0 otherwise.

$z_{jj'}^i$ : Binary variable that takes value 1 if the connection $(j, j')$ is used for container $i \in I$ along its route, and takes value 0 otherwise.

$n_i^k$ : Integer variable that indicates the number of $k$-blocks ($k \in K$) packed in container $i \in I$.

Note that, with these variables, the waste volume of those shipments that visit location $j$ with the purpose of collecting waste from modality $k$ can be expressed by means of $\sum_{i \in I} c \cdot n_i^k \cdot x_{ij}^k$

The objective function can be treated in three different ways depending on the problem.

1. When solving the classical BP, the objective is to minimize the number of containers used.

$$\text{Min } Z_1 := \sum_{i \in I} \sum_{k \in K} y_i^k$$

2. When solving the classical VRP, the objective is to minimize the total distance travelled.

$$\text{Min } Z_2 := \sum_{i \in I} \sum_{(j,j') \in A(J)} d_{jj'} \cdot z_{jj'}^i$$

3. We propose to apply a convex combination of both objectives with a parametric coefficient $\lambda \in (0,1)$ to be calibrated by municipal garbage collection services.

$$\text{Min } Z_3 := (1-\lambda) \cdot \sum_{i \in I} \sum_{k \in K} y_i^k + \lambda \cdot \sum_{i \in I} \sum_{(j,j') \in A(J)} d_{jj'} \cdot z_{jj'}^i$$

The nature of the variables used in the model allows us to build adequate programs to face different objectives. In our case, the following integer programming model inspired in the BP and VR schemes determines the deployment of routes for mobile eco-points for the selective collection of urban solid waste.

$$\text{Min } Z_3 := (1-\lambda) \cdot \sum_{i \in I} \sum_{k \in K} y_i^k + \lambda \cdot \sum_{i \in I} \sum_{(j,j') \in A(J)} d_{jj'} \cdot z_{jj'}^i$$

Subject to

$$\sum_{k \in K} y_i^k \geq 1 \qquad i \in I \tag{1}$$

$$\sum_{i \in I} \sum_{k \in K} x_{ij}^k \geq 1 \qquad j \in J \, (j \neq O) \tag{2}$$

$$y_i^k \leq n_i^k \leq |L| \cdot y_i^k \qquad k \in K; \; i \in I \tag{3}$$

$$\sum_{k \in K} n_i^k \leq |L| \qquad i \in I \qquad (4)$$

$$\sum_{j \in J} w_j^k \, x_{ij}^k \leq c \cdot n_i^k \qquad k \in K; \quad i \in I \qquad (5)$$

$$\sum_{i \in I} c \cdot n_i^k \cdot x_{ij}^k \geq w_j^k \qquad k \in K; \quad j \in J \, (j \neq O) \qquad (6)$$

$$\sum_{j' \in J / j' \neq O} z_{Oj'}^i = 1 \qquad i \in I \qquad (7)$$

$$\sum_{j \in J / j \neq O} z_{jO}^i = 1 \qquad i \in I \qquad (8)$$

$$\sum_{j' \in J /(j,j') \in A(J)} z_{jj'}^i - \sum_{j' \in J /(j',j) \in A(J)} z_{j'j}^i = 0 \qquad j \in J \, (j \neq O); \quad i \in I \qquad (9)$$

$$y_i^k \geq x_{ij}^k \, ; \quad \sum_{j' \in J /(j',j) \in A(J)} z_{j'j}^i \geq x_{ij}^k \qquad k \in K; \, j \in J \, (j \neq O); \, i \in I \qquad (10)$$

$$\sum_{(j,j') \in A(S)} z_{jj'}^i \leq |S| - 1 \qquad \{S : S \subseteq J, O \notin S, |S| \geq 2\}; \quad i \in I \qquad (11)$$

$$x_{ij}^k, \, y_i^k, \, z_{jj'}^i \in \{0,1\}; \, n_i^k \in \{1,2,...,|L|\}. \qquad (12)$$

Constraints (1) ensure that every route must collect at least some kind of waste. Constraints (2) establish that each location must be visited at least once. Constraints (3) - (4) guarantee consistency and an upper limit of the number of waste blocks of modality $k$ within each container. Constraints (5) ensure that the capacity of collecting type-$k$ waste in shipment $i$ is sufficient to satisfy the demand generated at nodes $j$ that are visited. Constraints (6) establish that the collection capacity of the total shipments that pass through point $j$ is enough to collect waste of each modality generated at that location. Constraints (7) - (9) guarantee the flow conservation at nodes, it is the classical constraint of VRP. Constraints (10) connect the decision variables used in the BP and VR blocks of the model. Constraints (11) are subtour elimination constraints. Constraints (12) indicate the nature of the variables used in the model.

## 4 Algorithm for solving the model

In the BP problem, items of different volumes must be packed into a finite number of bins (or containers), each of a fixed given volume, in order to minimize the number of bins used. The VRP addresses the determination of the optimal set of routes for a fleet of vehicles, in order to serve a given set of customers. Both models are of combinatorial nature and, in computational complexity theory, are classified as NP-hard problems. This fact justifies the use of algorithms that provide a good heuristic solution for the combined model. Furthermore, our proposed optimization model is non-linear, as can be seen in constraints (6). Therefore, in order to solve this complex problem, we propose an algorithm that simultaneously configures the containers and designs the routes that provide a good solution to our original problem.

When analyzing the set of restrictions in the model, we can appreciate the existence of quasi-separability between them, both in relation to the variables involved and the purpose pursued. Blocks (1) - (6) are aimed to established the configuration of the multiple bin container, while blocks (7) - (9) determine the most appropriate route to establish taking into account the existing demand for waste to be collected according to their types. As was previously pointed out, constraints (10) connect the decision variables used in the BP and VR blocks of the model.

If we algebraically manipulate the constraints (6), we can express them as follows:

$$\sum_{k \in K} \left( \sum_{i \in I} n_i^k \cdot x_{ij}^k \right) \geq \sum_{k \in K} \frac{w_j^k}{c}$$

Note that the following quotient determine the number of bin blocks required to collect the generation of $k$-waste at location $j$:

$$q_j^k = \left\lceil \frac{w_j^k}{c} \right\rceil$$

A matrix of $|J| \cdot |K|$ elements, whose individual components are the coefficients $q_j^k$ can be calculated according to the input data set. Note that:

- If $q_j^k > 0$ then node $j$ must be visited at least once for collecting the $k$-waste generated at location $j$. More specifically,
- If $0 < q_j^k \leq |L|$ then node $j$ does not need to be visited more than one time in order to collect all the $k$-waste generated at location $j$. The number of bins needed to collect all the waste of modality k located at point j could be concentrated in a single shipment i*, thus adapting the configuration of the container to the characteristics of the point to be visited (Adapted configuration strategy, option 2), or alternatively, it could be divided into several shipments that would repeatedly decrease the amount of waste to be collected. Among the multiple possible options for configuring the containers that would carry out these shipments is the one in which each type of waste is represented with a single bin in the container configuration (Fixed configuration strategy, option 1).

- If $q_j^k > |L|$ then node $j$ must necessarily be visited more than one time in order to collect all the $k$-waste generated at location $j$. In fact, quotient $\left\lceil \dfrac{q_j^k}{|L|} \right\rceil$ indicates the minimum number of required visits to carry out at location $j$. Since the total waste generated at location $j$ must be collected, according to constraints (6), we must assume that if the demand of collecting the $k$-waste generated at location $j$ is not satisfied by the visit of a first container, then a sequence of iterative visits must be implemented.

In the first phase of our algorithm, we will cover all the nodes $j$ that have the entire universe of specific residues. Hence, a container with a totally diversified configuration (that is, each block collects a different modality of residues from the rest of the blocks that configure the container) will travel from the depot node to node $j$ and will return following the shortest path.

In the second phase of our algorithm, there are only nodes where some types of waste are missing. Now we can choose between two strategies when configuring a multi-block container: to use identical containers, in which there are no repeated blocks in their configuration, or to use heterogeneous containers, where some block modalities have a greater presence, at the expense of others, when the configuration of the container is decided. Therefore, the efficiency of the algorithm that solves the problem depends on a correct choice in the order of action of two simple strategies:

**Option 1** (see Figure 3): To remove the plurality of blocks by using a single configuration for all collecting containers and extend their routes by visiting other nodes, until the bin blocks are filled with the corresponding specific materials.

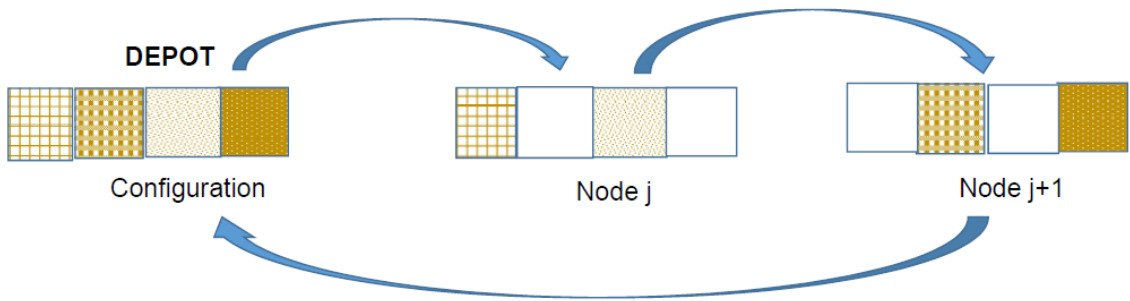

**Fig. 3**. Route for visiting nodes following the option 1 as configuration strategy.

**Option 2** (see Figure 4): To adapt the waste bin configuration of the visiting container to the characteristics of the node, since there is not a total variety of waste modalities at that point.

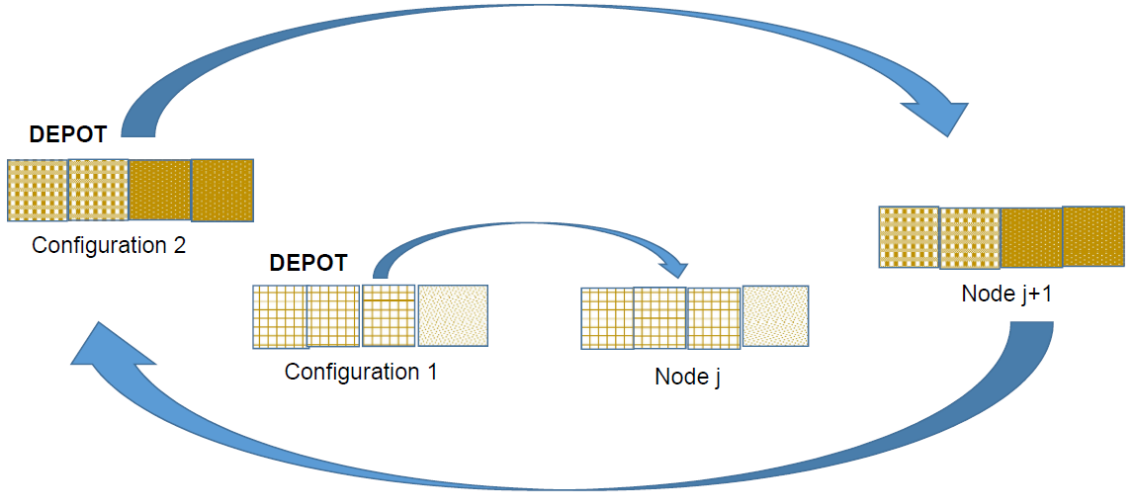

**Fig.4**. Two routes for visiting nodes with different container configurations adapted to the existing demand (option 2).

In the following heuristic both strategies are present.

- If only the main program is used, ignoring the call to the subroutine in step 4.1, we would actually be applying option 1, in which the container configuration is homogeneous for all vehicles. This situation would be applicable when the decision maker only has partial information on the demand at each point; that is, he/she knows that there is waste to be collected, but does not know its distribution in modalities, so optimizing the bin packaging is not a prerequisite.

- If we alternatively force the use of the subroutine in step 4.1, we would obtain a configuration of the multi-block container adapted to the existing demand in the nodes to be visited. Applying this second strategy (option 2) we would combine the search for shorter vehicle routes with the design of the most suitable packing of blocks for the container. If the decision maker had whole information on the demand at each point (that is, the distribution of the types of waste to be collected), he/she could optimize the container packaging before starting the route of the collecting container.

**HEURISTIC_1 (main program)**

1. **Sort** node set *J* according to the shortest distance from depot *O*.
2. **While**, at each node *j*, there exists all the universe of specific waste:

    **2.1 Generate** a route for a new visiting multi-block container *i* that follows a shortest path from depot *O* to node *j* along the street network.

    **2.2 Decrease** units from the unsatisfied demand at the visited node *j*.

3. **Identify** those nodes where some specific type of waste remains pending of collecting, but not for all types. **Let** *T* be this node set.
4. **While** there exist nodes *j* in *T* do

    **4.1 Select** the most appropriate configuration for a new multi-block bin *i* (by using subroutine HEURISTIC_2).

    **4.2 Generate** a route following a shortest path from depot *O* to node *j* along the street network.

    **4.3 Decrease** units of unsatisfied demand at each node *j*, according to the characteristics of the visiting multi-block bin.

    **4.4 If** the global demand of collecting specific waste has been satisfied, **Remove** node *j* from set *T*.

5. **End**.

**HEURISTIC_2 (subroutine to decide configuration of container *i* )**

**For** $k = 1$ to $|K|$ do

1. **If** $q_j^k \geq |L|$ **then** node *j* must be visited by means of a container *i* such that $n_i^k = |L|;\ \forall k \in K.$ **Return**

2. **If** $\sum_{k \in K} q_j^k \geq |L|$ **then** node *j* must be visited by means of a container *i* such that $\sum_{k \in K} n_i^k = |L|$. **Return**.

3. **If** $\sum_{k \in K} q_j^k < |L|$ **then** node *j* must be visited by means of a container *i* such that $n_i^k = q_j^k; \forall k \in K$. **Return**.

The tests carried out in the following section with randomly generated data have shown that option 2 involves the generation of shorter collecting routes, maintaining the number of routes required to guarantee the collection of all waste. Therefore, the use of option 2 is preferable, whenever possible.

## 5 Results and discussion

In order to validate the proposed optimization model, we have considered the Sioux Falls graph for different instances, varying the number of nodes and the quantity of waste produced at each node. In Figure 5, the Sioux Falls network with 24 nodes and 38 edges (76 directed arcs) is shown, where one depot is located at node 1 and the waste generating points are the points labeled 2-18. In Table 1 the results obtained when the problem of optimally deploying waste collection routes is implemented for several instances are shown. The objective function considered is the global distance traveled along the total of routes required to collect the totality of randomly generated waste.

We observe that for small size problems the solver achieves the optimal solution, improving the solutions obtained through Heuristic 1 (fixed configuration for all containers) and 2 (container configuration adapted to demand). However, when the network size increases the problem cannot be solved because of a problem of memory. As the experiments carried out have shown, when the exact model (solved by means of the GuRoBi 9.1.2 in Python solver and limiting the execution time to 2 hours) fails, both heuristic strategies achieve to find good solutions to the problem in a very short time. In particular, strategy 2, which combines the BP and VR optimization models, proves to be more efficient, as can be seen in Table 1.

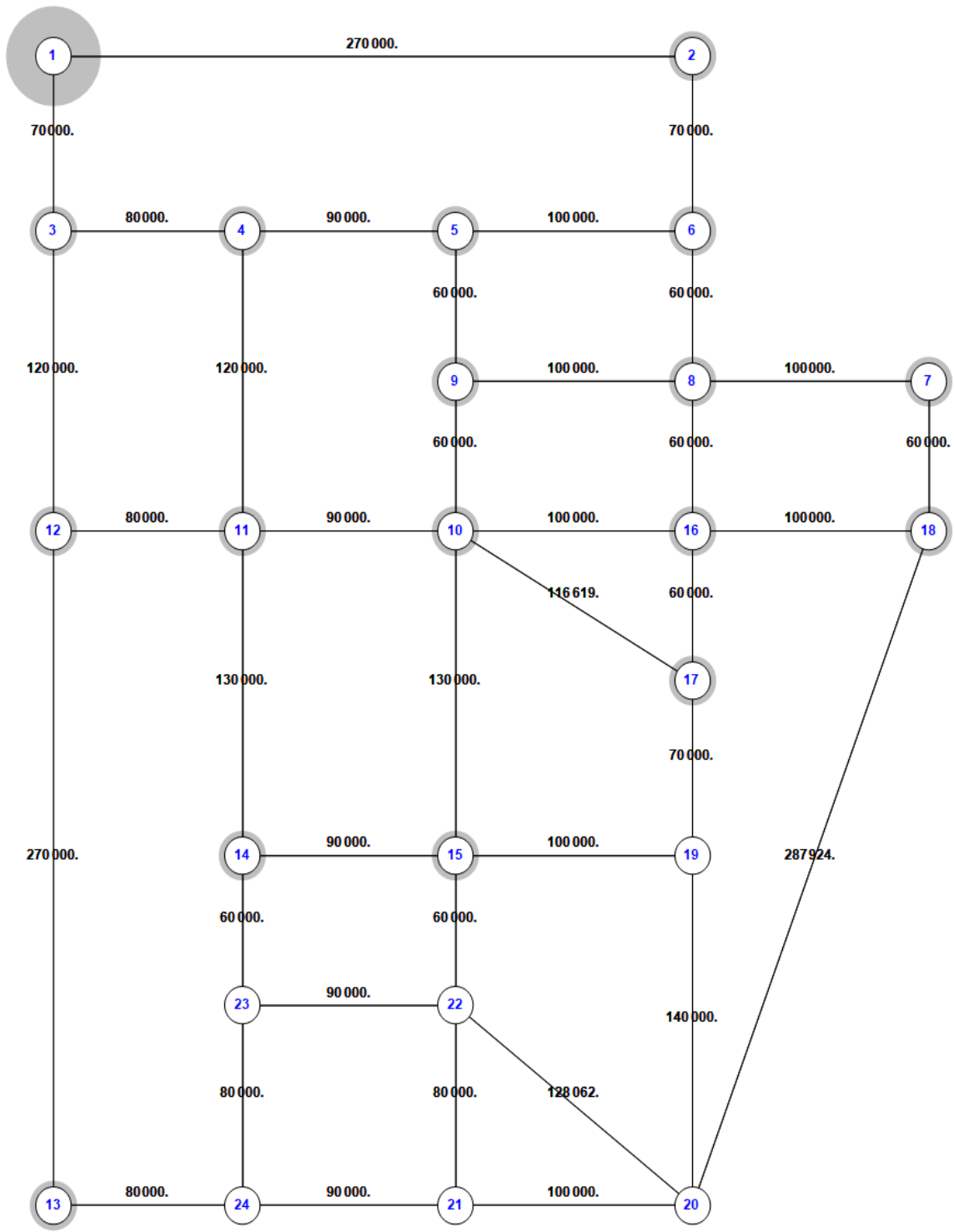

**Fig.5**. The Sioux Falls network with 24 nodes and 38 edges (76 directed arcs).

**Table 1:** Results obtained from the computational experience.

| Number of nodes | Objective function Exact model | CPU (seconds) Exact model | Heuristic 1 Fixed configuration | % | Heuristic 2 Adapted configuration | % |
|---|---|---|---|---|---|---|
| 10 | $7.020*10^6$ | 3.45 | $9.740*10^6$ | 38,75 | $7.180*10^6$ | 2,28 |
| 11 | $1.076*10^7$ | 4.82 | $1.282*10^7$ | 30,02 | $1.016*10^7$ | 3,04 |
| 12 | $9.260*10^6$ | 4.73 | $1.514*10^7$ | 63,50 | $9.680*10^6$ | 4,54 |
| 13 | $1.080*10^7$ | 5.90 | $1.552*10^7$ | 43,70 | $1.114*10^7$ | 3,15 |
| 14 | $1.240*10^7$ | 19.77 | $1.924*10^7$ | 55,16 | $1.290*10^7$ | 4,03 |
| 15 | $1.338*10^7$ | 76.79 | $1.974*10^7$ | 48,42 | $1.352*10^7$ | 1,05 |
| 16 | - | 7200 (3.08% Gap) | $1.986*10^7$ | - | $1.692*10^7$ | - |
| 18 | - | Out of memory | $2.294*10^7$ | - | $1.907*10^7$ | - |

In order to compare the efficiency of both heuristic strategies developed a computational experience has been carried out on a laboratory scenario composed of 39 nodes (sites identified as elements of node set *J*). Location of nodes along the street network is illustrated in Figure 6.

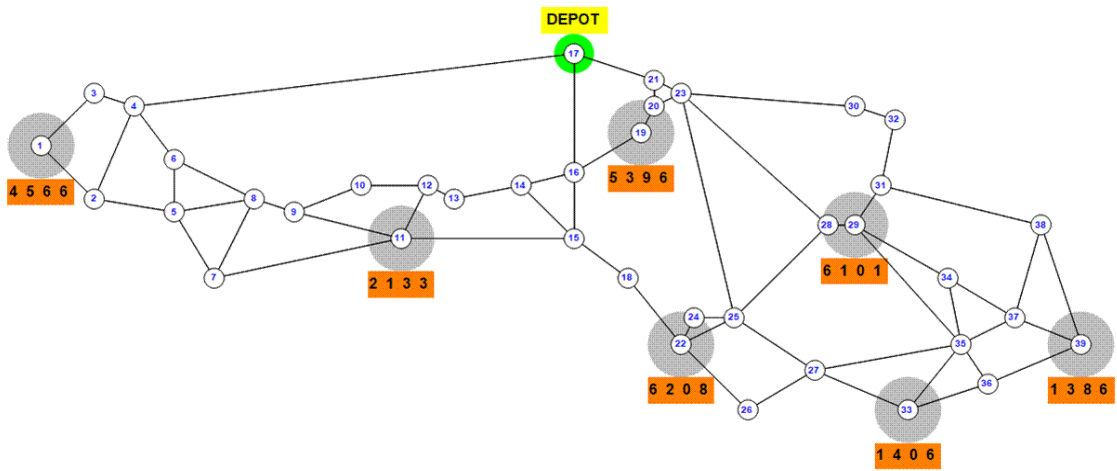

**Fig.6**. Network with 39 nodes and 7 waste collecting points.

Among the 39 nodes, 7 points have been selected (nodes labeled 1, 11, 19, 22, 29, 33 and 39) where the demand for a selective collection of waste will be located. The node labeled with number 17 represents the depot from where the routes of the mobile eco-points will start. Urban waste to be collected belongs to 4 categories and the quantity generated at each demand node is, for each category, a random number of integer nature belonging to the interval [1 kg, 9 kg]. The vehicle's collection capacity is limited to 4 kg in total, and its configuration may follow a homogeneous type (that is, 4 different blocks, each with the capacity to collect 1 kg) or an adapted type, being able to group the blocks in order to adapt it to the characteristics of the demand.

The first option is to use vehicles with the homogeneous configuration [1,1,1,1] in order to visit nodes where the demand has one representation of each modality. If in some node there were no units of a type, the vehicle should prolong its route to visit other nodes and, hence, complete its collection capacity before returning to the depot. In Figure 7 are illustrated the results in this case. Node 1 has the demand [2,0,1,1]; after this node is visited by the vehicle of homogeneous configuration, the vehicle would store the distribution [1,0,1,1] and the subsequent demand that would remain at the node would be [1,0,0,0].

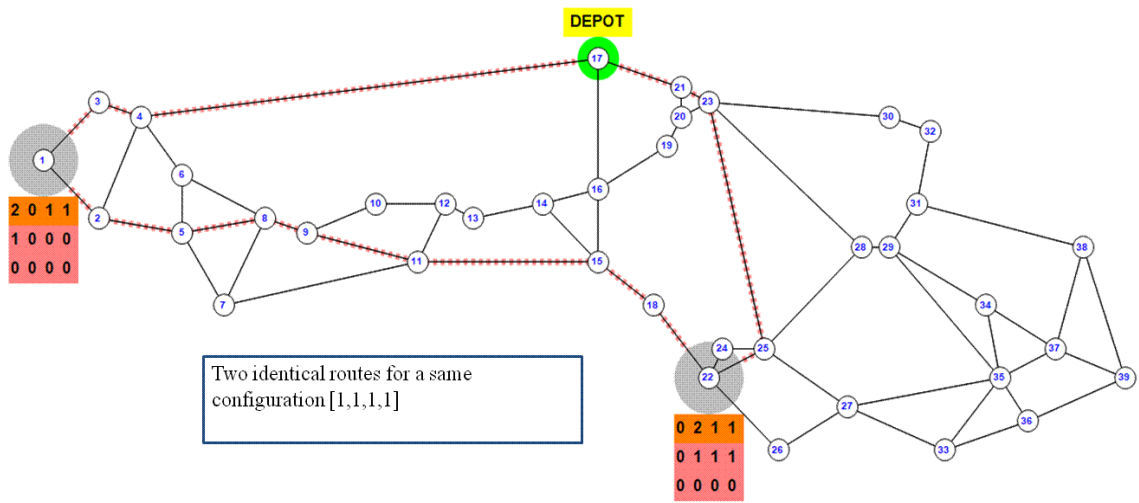

**Fig. 7**. Routes for the homogeneous configuration [1,1,1,1].

Therefore, this vehicle could continue its route in order to complete its loading capacity before returning to the depot; for example, visiting node 22 whose demand distribution is [0,2,1,1]. After visiting node 22, the vehicle would store [1,1,1,1] and then could return to the depot, whereas the demand of node 22 would now be [0,1,1,1].

A new vehicle of homogeneous configuration [1,1,1,1] could then sequentially revisit nodes 1 and 22, culminating the satisfaction of total demand with a second route. Applying this strategy, the cost of waste collection would be proportional to the total distance traveled by vehicles on both routes.

Alternatively, another strategy could be applied, where the distribution of the 4 blocks that the vehicle can transport was configured prior to starting the route from the depot. The chosen configuration would be determined by the existing demand at the nodes planned to visit. In the previously described case, where node 1 had the demand [2,0,1,1], it would be possible to dispatch a vehicle with the same configuration, since the sum of blocks (2+0+1+1) is exactly 4. In this way, the vehicle could travel directly from the depot to node 1 on the outward and return journey, without deviating from the shortest route.

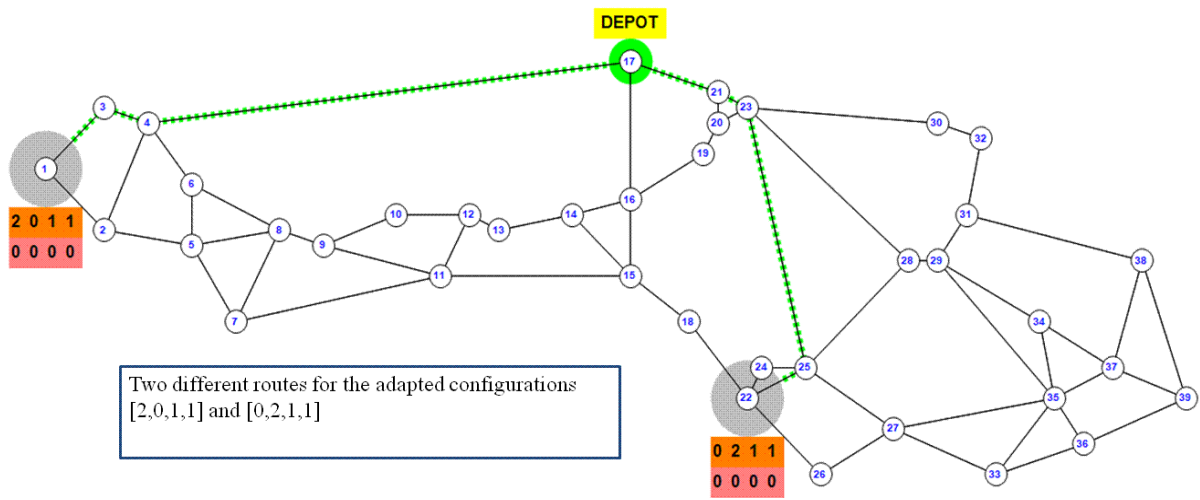

**Fig. 8**. Routes for adapted configurations of containers.

In Figure 8 it can be seen that after the vehicle visit, the demand located at node 1 is canceled ([0,0,0,0]). Analogously it would occur with node 22 whose initial demand was [0,2,1,1]. Therefore, the application of this second strategy allows apparently in a greater number of times the shortest routes for the travel of the mobile eco-points between the depot node and the demand points.

A computational experience, consisting of performing 30 experiments on the network in Figure 6 has been carried out in the above indicated terms in order to be able to compare the effectiveness of both strategies, which have been called respectively fixed (option 1) and adapted (option 2). Table 2 shows the results obtained.

**Table 2:** Results obtained from the computational experience.

| | Fixed configuration strategy | | Adapted configuration strategy | | Absolute improvements | | Relative Improvements | |
|---|---|---|---|---|---|---|---|---|
| n | # Routes | # Km | # Routes | # Km | # Routes | # Km | # Routes | # Km |
| 1 | 37 | 22587 | 37 | 19910 | 0 | 2677 | 0.00% | 11.85% |
| 2 | 32 | 24989 | 32 | 17505 | 0 | 7484 | 0.00% | 29.95% |
| 3 | 26 | 17959 | 27 | 15408 | -1 | 2551 | **-3.85%** | 14.20% |
| 4 | 29 | 18282 | 29 | 16743 | 0 | 1539 | 0.00% | 8.42% |
| 5 | 31 | 23388 | 31 | 18272 | 0 | 5116 | 0.00% | 21.87% |
| 6 | 34 | 25189 | 34 | 19476 | 0 | 5713 | 0.00% | 22.68% |
| 7 | 31 | 22758 | 30 | 17646 | 1 | 5112 | **3.23%** | 22.46% |
| 8 | 37 | 24855 | 36 | 19972 | 1 | 4883 | **2.70%** | 19.65% |
| 9 | 39 | 23740 | 39 | 21250 | 0 | 2490 | 0.00% | 10.49% |
| 10 | 29 | 20620 | 30 | 17749 | -1 | 2871 | **-3.45%** | 13.92% |
| 11 | 35 | 22384 | 35 | 19251 | 0 | 3133 | 0.00% | 14.00% |
| 12 | 35 | 23778 | 34 | 18737 | 1 | 5041 | **2.86%** | 21.20% |
| 13 | 33 | 23559 | 33 | 16722 | 0 | 6837 | 0.00% | 29.02% |
| 14 | 38 | 26480 | 38 | 22413 | 0 | 4067 | 0.00% | 15.36% |
| 15 | 27 | 16832 | 27 | 15083 | 0 | 1749 | 0.00% | 10.39% |
| 16 | 36 | 27400 | 36 | 19980 | 0 | 7420 | 0.00% | 27.08% |
| 17 | 33 | 24442 | 33 | 19117 | 0 | 5325 | 0.00% | 21.79% |
| 18 | 35 | 21937 | 34 | 17108 | 1 | 4829 | **2.86%** | 22.01% |
| 19 | 28 | 19032 | 28 | 15556 | 0 | 3476 | 0.00% | 18.26% |
| 20 | 32 | 23994 | 31 | 18744 | 1 | 5250 | **3.13%** | 21.88% |

| 21 | 29 | 21940 | 29 | 17738 | 0 | 4202 | 0.00% | 19.15% |
|---|---|---|---|---|---|---|---|---|
| 22 | 25 | 17548 | 25 | 14478 | 0 | 3070 | 0.00% | 17.49% |
| 23 | 31 | 21119 | 31 | 18045 | 0 | 3074 | 0.00% | 14.56% |
| 24 | 33 | 22022 | 32 | 16820 | 1 | 5202 | **3.03%** | 23.62% |
| 25 | 29 | 23679 | 28 | 15818 | 1 | 7861 | **3.45%** | 33.20% |
| 26 | 35 | 22148 | 34 | 18306 | 1 | 3842 | **2.86%** | 17.35% |
| 27 | 36 | 23894 | 36 | 19265 | 0 | 4629 | 0.00% | 19.37% |
| 28 | 24 | 17574 | 24 | 11680 | 0 | 5894 | 0.00% | 33.54% |
| 29 | 25 | 16908 | 25 | 13564 | 0 | 3344 | 0.00% | 19.78% |
| 30 | 35 | 22835 | 36 | 20814 | -1 | 2021 | **-2.86%** | 8.85% |

The columns of *Fixed configuration strategy* show the numbers of routes and total kilometers that a vehicle must travel to cover all the demand with this strategy (option 1). The columns of *Adapted configuration strategy* show the same with the adapted strategy (option 2). Finally, the columns of *Absolute improvements* and *Relative improvements* show the difference of the adapted strategy against the fixed strategy in terms of absolute cost and relative cost.

In all the experiments carried out, the adapted configuration strategy has improved the total cost invested in the determination of routes with respect to the results obtained by using the fixed configuration strategy. The improvement, in terms of percentage, is above 19% on average. In the number of vehicle routes dispatched from the depot, the improvement is inconclusive. As can be observed at Table 2, the number of routes is similar, regardless of the strategy used.

## 6 Conclusion

A methodology for the deployment of mobile multi-block containers for selective collection of urban solid waste has been proposed in this paper. The mathematical optimization model formulated for this purpose has been identified as a combined version of BP problem and the VR problem, whose computational complexity motivates the use of heuristics to face large real-life scenarios. Following that recommendation, a greedy algorithm has been developed to solve the proposed mathematical programming model. Two strategies have been identified for designing the configurations of the mobile multi-block containers that will visit the demand nodes. Sioux Falls network has been applied to show the state-of-art solvers are not capable to solve medium-size instances, but both heuristics provide good solutions for small-size instances.

To decide the most efficient strategy for implementing a solving algorithm, a computational experience has been carried out on a laboratory instance. Results show that to adapt the configuration of the mobile multi-block container to the characteristics of the node to be visited must be a priority.

This work can provide a useful information to environment engineers and operators in the WM field, like possible performance enhancement modifications to urban waste collection.

## ACKNOWLEDGEMENTS

This work was in part supported by the Ministerio de Economía y Competitividad (Spain)/FEDER under grant MTM2015-67706-P, and by the Ministerio de Investigación (Spain)/FEDER under grant PID2019-106205GB-I00. This support is gratefully acknowledged.